\documentclass[11pt,a4paper]{article}
\usepackage{epsfig}
\usepackage{epic}
\usepackage{eepic}
\usepackage{bm}
\usepackage{amssymb}
\usepackage{amsmath}
\usepackage{caption}
\usepackage{subcaption}
\usepackage{colortbl}
\usepackage{color}
\usepackage{hyperref}

\newtheorem{theorem}{Theorem}[section]

\newtheorem{lemma}[theorem]{Lemma}

\renewcommand{\Re}{\mathrm {I\!R}}
\newcommand{\naturalNumbers}{\mathrm {I\!N}}

\title{Regularization of ill-posed point neuron models}
\author{Bj{\o}rn Fredrik Nielsen\thanks{Faculty of Science and Technology, Norwegian University of Life Sciences, P.O. Box 5003, NO-1432 {\AA}s, NORWAY. Email: bjorn.f.nielsen@nmbu.no.}
}

\begin{document}
\maketitle

\begin{abstract}
Point neuron models with a Heaviside firing rate function can be ill-posed. That is, the {\em initial-condition}-to-{\em solution} map might become discontinuous in finite time.
If a Lipschitz continuous, but steep, firing rate function is employed, then standard ODE theory implies that such models are well-posed and can thus, approximately, be solved with finite precision arithmetic. We investigate whether the solution of this well-posed model converges to a solution of the ill-posed limit problem as the steepness parameter, of the firing rate function, tends to infinity. Our argument employs the Arzel\`{a}-Ascoli theorem and also yields the existence of a solution of the limit problem. However, we only obtain convergence of a subsequence of the regularized solutions. This is consistent with the fact that we show that models with a Heaviside firing rate function can have several solutions. Our analysis assumes that the Lebesgue measure of the time the limit function, provided by the Arzel\`{a}-Ascoli theorem, equals the threshold value for firing, is zero. If this assumption does not hold, we argue that the regularized solutions may not converge to a solution of the limit problem with a Heaviside firing function.  
%
\end{abstract}


\section{\label{Introduction}Introduction}
In this paper we analyze some mathematical properties of the following classical point neuron model: 
\begin{align}
\label{model-eq1}
& \tau_{i} u'_i(t) = -u_i(t) + \sum_{j=1}^{N} \omega_{i,j} S_{\beta}[u_j(t)-u_{\theta}] + q_i(t), \, t \in (0,T],   \\
\label{model-eq2}
& u_i(0)=u_{\mathrm{init},i}, 
\end{align}
for $i=1,2,\ldots,N$, where 
\begin{align}
&u_i(t) \in \Re, \, t \in [0,T], \, i=1,2,\ldots,N, \nonumber \\
&q_i(t) \in \Re, \, t \in (0,T], \, i=1,2,\ldots,N, \nonumber \\
&u_{\mathrm{init},i} \in \Re, \, i=1,2,\ldots,N, \nonumber \\
&u_{\theta} \in \Re,  \nonumber \\
&\omega_{i,j} \in \Re, \, i,j=1,2,\ldots,N, \nonumber \\
&\tau_i \in \Re_+, \, i=1,2,\ldots,N, \nonumber\\
&\beta=1,2,\ldots,\infty, \nonumber \\
&S_{\beta}[x]  \mbox{ is an approximation of the Heaviside function } H[x], \nonumber\\
&S_{\infty}[x]=H[x]. \nonumber
\end{align}
Here, $u_{i}(t)$ represents the unknown electrical potential of the $i$th unit in a network of $N$ units. The nonlinear function $S_{\beta}$ is called the {\em firing rate function}, $\beta$ is the {\em steepness parameter} of $S_{\beta}$, $u_{\theta}$ is the {\em threshold value for firing}, $\{ \omega_{ij} \}$ are the connectivities, $\{ \tau_i \}$ are membrane time constants and $\{ q_i(t) \}$ model the external drive/external sources, see, e.g., \cite{B2012,E1998,FVG2009} for further details. 

In computational neuroscience one often employs a steep sigmoid, or Heaviside, firing rate function $S_{\beta}$. This is due to both electrophysiological properties and mathematical convenience\footnote{Amari \cite{Amari} analyzed the stationary solutions of neural field equations when $\beta = \infty$. }. 
Unfortunately, the {\em initial-condition}-to-{\em solution} map for (\ref{model-eq1})-(\ref{model-eq2}) can become discontinuous, in finite time, if a Heaviside firing rate function is used \cite{NW2016}. Such models are thus virtually impossible to solve with finite precision arithmetic \cite{EHN1996,Wikipedia}. Also, in the steep, but Lipschitz continuous, firing rate regime, the error amplification
can be extreme, even though a minor perturbation of the initial condition does not change which neurons that fire. It is important to note that this ill-posed nature of the model is a fundamentally different mathematical property than the possible existence of unstable equilibria, which typically also occur if a firing rate function with moderate steepness is used, see \cite{NW2016} for further details. 

The solution of (\ref{model-eq1})-(\ref{model-eq2}) depends on the steepness parameter $\beta$. That is, 
\[
u_i(t) = u_{\beta,i}(t), \, i=1,2,\ldots,N,
\]
and the purpose of this paper is to analyze the limit process $\beta \rightarrow \infty$. 
This investigation is motivated by the fact that the stable numerical solution of an ill-posed problem is very difficult, if not to say impossible, see, e.g., \cite{EHN1996,Wikipedia}. Consequently, such models must be regularized to obtain a sequence of well-posed equations which, at least in principle, can be approximately solved by a computer. Also,  steep firing rate functions, or even the Heaviside function, are often used in simulations. It is thus necessary to explore whether the
limit process $\beta \rightarrow \infty$ is mathematically sound. Similar type of studies, using different techniques, are  presented in \cite{OPW2013b,OPKS2016} for the stationary solutions of neural field models. 

 
In sections \ref{uniformily-bounded-equicontinuous} and \ref{convergence-expected-limit} we use the Arzel\`{a}-Ascoli theorem to analyze the properties of the sequence $\{ \mathbf{u}_{\beta} \}$, where 
\begin{equation}
\label{u-beta}
\mathbf{u}_{\beta} (t)= \left( u_{\beta,1}(t), \, u_{\beta,2}(t), \, \ldots, \, u_{\beta,N}(t) \right)^T.
\end{equation}
More specifically, we prove that this sequence has at least one subsequence which converges uniformly to a limit 
\[
\mathbf{v}(t)=(v_1(t), \, v_2(t), \ldots, \, v_N(t))^T,
\] 
and that this limit satisfies the integral/Volterra version of (\ref{model-eq1})-(\ref{model-eq2}) with $S_{\beta}=S_{\infty}$, provided that the Lebesgue measure of the time one, or more, of the component functions of $\mathbf{v}$ equals the threshold value $u_{\theta}$ for firing, is zero. 
Furthermore, in section \ref{threshold-advanced-limits} we argue that, if $\mathbf{v}$ does not satisfy this threshold property, then this function will not necessarily solve the limit problem. 

According to the Picard-Lindel\"{o}f theorem \cite{WikipediaPicard}, (\ref{model-eq1})-(\ref{model-eq2}) has a unique solution, provided that 
$\beta < \infty$ and that the assumptions presented in the next section hold. In section \ref{uniqueness} we show that this uniqueness feature is \underline{not} necessary inherited by the limit problem obtained by employing a Heaviside firing rate function. It actually turns out that different subsequence of $\{ \mathbf{u}_{\beta} \}$ can converge to different solutions of (\ref{model-eq1})-(\ref{model-eq2}) with $S_{\beta}=S_{\infty}$. This is explained in section \ref{convergence-entire-sequence}, which also contains a result addressing the convergence of the entire sequence  $\{ \mathbf{u}_{\beta} \}$.  


For the sake of easy notation, we will sometimes write (\ref{model-eq1})-(\ref{model-eq2}) in the form 
\begin{align}
\label{point-neuron-model}
\mathbf{\tau} \mathbf{u}'(t) &= - \mathbf{u}(t) + \mathbf{\omega}  S_{\beta}[\mathbf{u}(t)-\mathbf{u_{\theta}}] + \mathbf{q}(t), \, t \in (0,T],\\
\label{initial-condition}
\mathbf{u}(0)&=\mathbf{u}_{\mathrm{init}},
\end{align}
where 
\begin{align}
&\mathbf{u}(t)=\mathbf{u}_{\beta} (t) \in \Re^N, \, t \in [0,T], \mbox{ see (\ref{u-beta})}, \nonumber \\
&\mathbf{q}(t)=\left( q_1(t), \, q_2(t), \, \ldots, \, q_N(t) \right)^T \in \Re^N, \, t \in (0,T], \nonumber \\
&\mathbf{u_{\theta}}= \left(u_{\theta}, \, u_{\theta}, \ldots, \, u_{\theta} \right)^T \in \Re^N, \nonumber \\
&\mathbf{u}_{\mathrm{init}} = \left(u_{\mathrm{init},1}, \, u_{\mathrm{init},2} , \ldots, \, u_{\mathrm{init},N} \right)^T \in \Re^N, \nonumber \\
&\mathbf{\omega} = [\omega_{i,j}] \in \Re^{N \times N}, \nonumber \\
&\mathbf{\tau} = \mathrm{diag}(\tau_1,\, \tau_2, \ldots, \tau_N) \in \Re^{N \times N} \mbox{ is diagonal}, \nonumber\\
&S_{\beta}[\mathbf{x}]=(S_{\beta}[x_1], \ldots,S_{\beta}[x_N])^T, \,
\mathbf{x} = \textcolor{black}{(x_1, \ldots, x_N)^T} \in \Re^N. \label{vector-notation1}
\end{align}

\section{Assumptions}
Throughout this text we use the standard notation
\begin{equation}
\label{norm}
\|\mathbf{x}\|_{\infty}=\max \limits_{1\leq i\leq N}|x_{i}|,\quad \mathbf{x}=(x_{1},...,x_{N})  \in \Re^N.
\end{equation}

Concerning the sequence $\{ S_{\beta} \}$ of finite steepness firing rate functions, we make the following assumption. 

\vspace*{0.2cm}
\noindent{\bf Assumption \cal{A}} \\
\emph{We assume that 
\begin{description}
\item[a)] $S_{\beta}$, $\beta \in \naturalNumbers$, is Lipschitz continuous, 
\item[b)]  $0 \leq S_{\beta}(x) \leq 1$, $x \in \Re, \, \beta \in \naturalNumbers$, 
\item[c)] for every pair of positive numbers $(\epsilon, \, \delta)$ there exists $Q \in \naturalNumbers$ such that 
\begin{align}
\label{H_app1}
& |S_{\beta}(x)| < \epsilon \quad \mbox{for } x < -\delta \mbox { and } \beta > Q, \\
\label{H_app2}
& |1-S_{\beta}(x)| < \epsilon \quad \mbox{for } x > \delta \mbox { and } \beta > Q.
\end{align} 
\end{description}
}
Reasonable/sound approximations of the Heaviside function satisfy {\bf \cal A}. 
For example, if $S_{\beta}$ is nondecreasing (for every $\beta \in \naturalNumbers$), a) and b) hold and $\{ S_{\beta} \}$ converges pointwise to the Heaviside function, then {\bf {\cal A}} holds. 
Also, if assumption {\bf {\cal A}} is satisfied and 
$\lim_{\beta \rightarrow \infty} S_{\beta}(0) = S_{\infty}(0)=H(0)$, then $\{ S_{\beta} \}$  converges pointwise to the Heaviside function. Many continuous sigmoid approximations of the Heaviside function obey {\bf \cal{A}}. For example, 
\begin{align}
&S(x) = \frac{1}{2}(1+\tanh(x)) \label{firingtanh1}, \\
&S_{\beta}[x] = S(\beta x). \label{firingtanh2}
\end{align}

We will consider a slightly more general version of the model than (\ref{point-neuron-model})-(\ref{initial-condition}). More specifically, we allow the source term to depend on the steepness parameter, $\mathbf{q}=\mathbf{q_{\beta}}$, but in such a way that the following assumption holds. 

\vspace*{0.2cm}
\noindent{\bf Assumption \cal{B}} \\
\emph{We assume that $\mathbf{q}_{\beta} (t)$, $t \in [0,T]$, $\beta \in \naturalNumbers \cup \{ \infty \}$ is continuous and that 
\begin{align}
\label{q-bounded}
& \sup_{\beta \in \naturalNumbers, \, t \in [0,T]} \| \mathbf{q_{\beta}} (t) \|_{\infty}  \leq B < \infty, \, B \in \Re,  \\
\label{limit-q}
& \lim_{\beta \rightarrow \infty} \mathbf{q_{\beta}} (t) = \mathbf{q_{\infty}} (t), \, t \in [0,T], \\
\label{limit-int-q}
&  \lim_{\beta \rightarrow \infty} \int_0^t \mathbf{q_{\beta}}(s) \, ds = \int_0^t \mathbf{q_{\infty}}(s) \, ds, 
 \, t \in [0,T].
\end{align}
}
Allowing the external drive to depend on the steepness parameter, makes it easier to construct illuminating examples. 
Please note that our theorems also will hold for the simplest case, i.e. when $\mathbf{q}$ does not change as $\beta$
increases. 

%

\noindent In this paper we will assume that assumptions {\bf \cal{A}} and {\bf \cal{B}} are satisfied. 

\section{Uniformly bounded and equicontinuous}
\label{uniformily-bounded-equicontinuous}
In order to apply the Arzel\`{a}-Ascoli theorem we must show that $\{ \mathbf{u}_{\beta} \}$ constitute a family of uniformly bounded and equicontinuous functions. (For the sake of simple notation, we will write $u_i$ and $q_i$, instead of $u_{\beta,i}$ and $q_{\beta,i}$, for the component functions of $\mathbf{u}_{\beta}$ and $\mathbf{q}_{\beta}$, respectively). Multiplying
\[
u'_i(s) + \tau_{i}^{-1} u_i(s) = \tau_{i}^{-1} \sum_{j=1}^{N} \omega_{i,j} S_{\beta}[u_j(s)-u_{\theta}] + \tau_{i}^{-1} q_i(s)
\]
with $e^{\tau_{i}^{-1} s}$ yields that
\[
\left[ u_i(s) e^{\tau_{i}^{-1} s} \right]' =  e^{\tau_{i}^{-1} s} \tau_{i}^{-1} \sum_{j=1}^{N} \omega_{i,j} S_{\beta}[u_j(s)-u_{\theta}] + e^{\tau_{i}^{-1} s} \tau_{i}^{-1} q_i(s)
\]
and by integrating 
\[
u_i(t) e^{\tau_{i}^{-1} t} = u_i(0) + \int_0^t e^{\tau_{i}^{-1} s} \tau_{i}^{-1} \sum_{j=1}^{N}   \omega_{i,j} S_{\beta}[u_j(s)-u_{\theta}] \, ds + \int_0^t e^{\tau_{i}^{-1} s} \tau_{i}^{-1} q_i(s) \, ds.
\]
Hence, since $S_{\beta}[x] \in [0,1]$ and we assume that $\tau_i > 0$ for $i=1,2,\ldots,N$,
\begin{align*}
|u_i(t)| e^{\tau_{i}^{-1} t} &\leq |u_i(0)| + \sum_{j=1}^{N}  |\omega_{i,j}| \int_0^t e^{\tau_{i}^{-1} s} \tau_{i}^{-1} \, ds + \sup_{s \in [0,T]} |q_i(s)| \int_0^t e^{\tau_{i}^{-1} s} \tau_{i}^{-1} \, ds \\
&= |u_i(0)| + \left( \sum_{j=1}^{N}  |\omega_{i,j}| + \sup_{s \in [0,T]} |q_i(s)| \right) \left( e^{\tau_{i}^{-1} t} - 1 \right) \\
&\leq |u_i(0)| + \left( \sum_{j=1}^{N}  |\omega_{i,j}| + B \right) \left( e^{\tau_{i}^{-1} t} - 1 \right), 
\, t \in (0,T], 
\end{align*}
where the last inequality follows from (\ref{q-bounded}).
This implies that
\begin{equation}
\label{uniformly-bounded}
\| \mathbf{u}_{\beta}(t) \|_{\infty}  \leq \| \mathbf{u}_{\mathrm{init}} \|_{\infty} + \max_i \left( \sum_{j=1}^{N}  |\omega_{i,j}|  \right) + B, \, t \in [0,T].
\end{equation}
Since the right-hand-side of (\ref{uniformly-bounded}) is independent of $\beta$ and $t$, we conclude that the sequence $\{ \mathbf{u}_{\beta} \}$ is uniformly bounded. Next, from the bound (\ref{uniformly-bounded}), the model equation (\ref{point-neuron-model}), assumption (\ref{q-bounded}) and the assumption that $S_{\beta}[x] \in [0,1]$ we find that also $\{ \mathbf{u}'_{\beta} \}$ is uniformly bounded. It therefore follows from the Mean Value theorem that $\{ \mathbf{u}_{\beta} \}$ is a set of equicontinuous functions.

The Arzel\`{a}-Ascoli theorem \cite{WikipediaArzelaAscoli, Griffel, Royden} now asserts that there is a uniformly convergent subsequence $\{ \mathbf{u}_{\beta_k} \}$:
\begin{equation}
\label{limit-function}
\mathbf{v} = \lim_{k \rightarrow \infty} \mathbf{u}_{\beta_k}.
\end{equation}
According to standard ODE theory, $\mathbf{u}_{\beta}$ is continuous for $\beta=1,2,\ldots < \infty$, and hence the uniform  convergence implies that also $\mathbf{v}$ is continuous.

\subsection{Threshold terminology}
\label{threshold-terminology}
As we will see in subsequent sections, whether we can prove that $\mathbf{v}$ actually solves the limit problem with a Heaviside firing rate function, depends on $\mathbf{v}$'s threshold properties. The following concepts turn out to be useful. 

For a  vector-valued function $\mathbf{z} = (z_1,z_2,\ldots,z_N)^T: [0,T] \rightarrow \Re^N$ we define
\begin{equation}
\label{ms}
m(s;\mathbf{z}) = \min_{j \in \{1,2,\ldots, N\}} |z_j(s)-u_{\theta}|, \, s \in [0,T].
\end{equation}

\vspace*{0.2cm}
\noindent
{\bf Definition} {\bf \em (Threshold simple)} \\
{\em A measurable vector-valued function $\mathbf{z}:[0,T] \rightarrow \Re^N$ is threshold simple if the Lebesgue measure of the set 
\begin{equation}
\label{zero-measure}
Z(\mathbf{z}) = \{s  \in [0,T] \, | \, m(s;\mathbf{z})=0 \}
\end{equation}
is zero, i.e. $|Z(\mathbf{z})|=0$. 
}

\vspace*{0.2cm}
\noindent
{\bf Definition} {\bf \em (Extra threshold simple)} \\
{\em A measurable vector-valued function $\mathbf{z}:[0,T] \rightarrow \Re^N$ is extra threshold simple if there exist \underline{open} intervals
\[
I_l = (a_l, a_{l+1}), \, l=1,2,\ldots, L,
\]
such that
\begin{eqnarray*}
&& a_1=0, \, a_{L+1} = T, \\
&& m(s;\mathbf{z}) \neq 0 \quad \forall s \in \bigcup_{l=1}^L I_l.
\end{eqnarray*}}

With words, $\mathbf{z}$ is extra threshold simple if $|Z(\mathbf{z})|=0$ and the component functions of $\mathbf{z}$ only attains the threshold value for firing $u_{\theta}$ a finite number of times during $[0,T]$. 
\vspace*{0.2cm}

\section{Convergence to the expected limit}
\label{convergence-expected-limit}
\subsection{Preparations}
We will prove that the limit $\mathbf{v}$ in (\ref{limit-function}) solves the integral form of (\ref{point-neuron-model})-(\ref{initial-condition}) with $S_{\infty} = H$, the Heaviside function, provided that $\mathbf{v}$ is threshold simple. 
The inhomogeneous nonlinear Volterra equation associated with (\ref{point-neuron-model})-(\ref{initial-condition}) reads: 
\begin{align}
\nonumber
\mathbf{\tau} \mathbf{u}_{\beta_k}(t) - \mathbf{\tau} \mathbf{u}_{\mathrm{init}}=& - \int_0^t \mathbf{u}_{\beta_k}(s) \, ds \\
\nonumber
& + \int_0^t \mathbf{\omega}  S_{\beta_k}[\mathbf{u}_{\beta_k}(s)-\mathbf{u_{\theta}}] \, ds \\
\label{integral-form}
& + \int_0^t \mathbf{q_{\beta_k}}(s) \, ds, \, t \in [0,T],
\end{align}
where we consider the equations satisfied by the subsequence $\{ \mathbf{u}_{\beta_k} \}$, see (\ref{limit-function}). We will analyze the convergence of the entire sequence in section \ref{convergence-entire-sequence}. Note that we use the notation 
\[
\int_0^t \mathbf{u}_{\beta_k}(s) \, ds = 
\left( 
\int_0^t u_{\beta_k,1}(s) \, ds, \, \int_0^t u_{\beta_k,2}(s) \, ds, \ldots, \, \int_0^t u_{\beta_k,N}(s) \, ds 
\right)^T
\]
etc. in (\ref{integral-form}), see also (\ref{u-beta}) and (\ref{vector-notation1}).

%

The uniform convergence of $\{ \mathbf{u}_{\beta_k} \}$ to $\mathbf{v}$ implies that the left-hand-side and the first term on the right-hand-side of (\ref{integral-form}) converge to the "expected" limits as $k \rightarrow \infty$. Also, due to assumption (\ref{limit-int-q}), the third term on the right-hand-side does not require any extra attention. We will thus focus on the second term on the right-hand-side of (\ref{integral-form}).

Let, for $t \in [0,T]$ and $\delta > 0$, 
\begin{align}
\label{p-def}
& p(\delta;t) = \left\{ s \in [0,t] \, | \, m(s;\mathbf{v}) > \delta \right\}, \\
\label{r-def}
& r(\delta;t) = [0,t] \setminus p(\delta;t),
\end{align}
where $m(s;\mathbf{v})$ is defined in (\ref{ms}), and $\mathbf{v}$ is the limit in (\ref{limit-function}). 
We note that, provided that $\delta >0$ is small, the set $r(\delta;t)$ contains the times where at least one of the components of $\mathbf{v}$ is close to the threshold value $u_{\theta}$ for firing. The following lemma turns out to be crucial for our analysis of 
the second term on the right-hand-side of (\ref{integral-form})
\begin{lemma}
If the limit function $\mathbf{v}$ in (\ref{limit-function}) is threshold simple, then 
\begin{equation}
\label{corollary-conjecture}
\lim_{\delta \rightarrow 0^+} |r(\delta;t)| = 0, \, t \in [0,T], 
\end{equation}
where $|r(\delta;t)|$ denotes the Lebesgue measure of the set $r(\delta;t)$. 
\end{lemma}
\subsubsection*{Proof}
\begin{itemize}
\item Since $\mathbf{v}$ is the uniform limit of a sequence of continuous functions, $\mathbf{v}$ is continuous and hence measurable.
\item If $\mathbf{v}$ is threshold simple, then 
\begin{equation}
\label{zero-measure2}
|Z(\mathbf{v})|=0,
\end{equation}
see (\ref{zero-measure}).
\item Let $t \in [0,T]$ be arbitrary. 
\item Assume that $\lim_{\delta \rightarrow 0^+} |r(\delta;t)| \neq 0$, or that this limit does not exist. 
\item Then $\exists \, \tilde{\epsilon} > 0$ such that there is a sequence $\{ \delta_n \}$ satisfying  
\begin{align*}
& 0 < \delta_{n+1} < \delta_n \quad \forall n \in \naturalNumbers, \\
& \lim_{n \rightarrow \infty} \delta_n = 0, \\
& |r(\delta_n;t)| > \tilde{\epsilon} \quad \forall n \in \naturalNumbers. 
\end{align*}
\item By construction, 
\[
r(\delta_1;t) \supset r(\delta_2;t) \supset \ldots \supset r(\delta_n;t) \supset \ldots, 
\]
and $|r(\delta_1;t)| \leq T < \infty$. Hence, 
\[
\left| \bigcap_{n=1}^{\infty} r(\delta_n;t) \right| = \lim_{n \rightarrow \infty} |r(\delta_n;t)| \geq \tilde{\epsilon} > 0, 
\]
see, e.g., \cite{Royden} (page 62). Since the sequence $\{ |r(\delta_n;t)| \}$ is nonincreasing and bounded below, $\lim_{n \rightarrow \infty} |r(\delta_n;t)|$ exists.  
\item Next, 
\[
s \in \bigcap_{n=1}^{\infty} r(\delta_n;t) \, \Rightarrow \, m(s;\mathbf{v}) \leq \delta_n \; \forall n \, \Rightarrow \, m(s;\mathbf{v})=0 \, \Rightarrow \, s \in Z(\mathbf{v}),  
\]  
i.e. 
\[
\bigcap_{n=1}^{\infty} r(\delta_n;t) \subset Z(\mathbf{v}). 
\]
\item Hence, 
\[
| Z(\mathbf{v}) | \geq \left| \bigcap_{n=1}^{\infty} r(\delta_n;t) \right| \geq \tilde{\epsilon} > 0, 
\]
which contradicts (\ref{zero-measure2}).
\end{itemize}
\rule{3mm}{3mm}

\subsection{Convergence of the integral}
\begin{lemma}
If the limit $\mathbf{v}$ in (\ref{limit-function}) is threshold simple, then 
\begin{equation}
\label{C5}
\lim_{k \rightarrow \infty}
\int_0^t \mathbf{\omega} S_{\beta_k}[\mathbf{u}_{\beta_k}(s)-\mathbf{u_{\theta}}] \, ds =
\int_0^t \mathbf{\omega} S_{\infty}[\mathbf{v}(s)-\mathbf{u_{\theta}}] \, ds, \, t \in [0,T].
\end{equation}
\end{lemma}
\subsubsection*{Proof}
Let $t \in [0,T]$ and $\tilde{\epsilon} > 0$ be arbitrary, and define
\[
C= \max_{i \in \{1,2,\ldots,N\}} \left( \sum_{j=1}^{N}  |\omega_{i,j}| \right).
\]
From (\ref{corollary-conjecture}) we know that there exists $\Delta > 0$ such that
\begin{equation}
\label{C0}
|r(2\delta;t)| < \frac{\tilde{\epsilon}}{2C}, \, 0 < \delta < \Delta.
\end{equation}
Choose a $\delta$ which satisfies $0 < \delta < \Delta$. 
According to assumption {\bf \cal{A}}, for this $\delta$ and 
\begin{equation}
\label{C0.01}
\epsilon=\frac{\tilde{\epsilon}}{2TC},
\end{equation}
there exists $Q \in \naturalNumbers$ 
such that (\ref{H_app1}) and (\ref{H_app2}) hold.  

Recall that $\beta_1, \, \beta_2, \ldots, \, \beta_k, \, \ldots$ are the values for the steepness parameter associated with the convergent subsequence $\{ \mathbf{u}_{\beta_k} \}$ in (\ref{limit-function}). 
Let $K \in \naturalNumbers$ be such that
\begin{align}
\label{C0.1}
&\beta_K > Q, \\
\label{C1}
&\sup_{s \in [0,T]} \| \mathbf{u}_{\beta_k}(s) - \mathbf{v}(s) \|_{\infty} < \delta, \, k > K.
\end{align}
The existence of such a $K$ is assured by the uniform convergence of $\{ \mathbf{u}_{\beta_k} \}$ to $\mathbf{v}$.
From the definition of the set $p(2 \delta;t)$, see (\ref{p-def}) and (\ref{ms}), 
\begin{equation}
\label{C2}
m(s; \mathbf{v}) = \min_{ j \in \{1,2,\ldots, N \}} |v_j(s)-u_{\theta}| > 2 \delta > \delta, \, s \in p(2 \delta;t),
\end{equation}
and from (\ref{C1}), and the triangle inequality, it follows that
\begin{equation}
\label{C3}
\min_{ j \in \{1,2,\ldots, N \}} |u_{\beta_k,j}(s)-u_{\theta}| > \delta, \, s \in p(2 \delta;t) \mbox{ and } k > K.
\end{equation}

From (\ref{C1}) and (\ref{C2}) we find that 
\[
(v_j(s)-u_{\theta}) \cdot (u_{\beta_k,j}(s)-u_{\theta}) > 0, \, s \in p(2 \delta;t), j \in \{1,2,\ldots, N \}, \, k > K.  
\]
Also, because of the properties of the Heaviside function,  
\[
S_{\infty}(v_j(s)-u_{\theta}) = 
\left\{
\begin{array}{ll}
1, & v_j(s)-u_{\theta} \geq \delta, \\
0 & v_j(s)-u_{\theta} \leq -\delta, 
\end{array}
\right.
\]
$j \in \{1,2,\ldots, N \}$. 
%
%
Consequently, due to (\ref{C0.1}) and assumption {\bf \cal{A}}, see (\ref{H_app1}) and (\ref{H_app2}), we find that
\[
|S_{\beta_k}[u_{\beta_k,j}(s)-u_{\theta}] - S_{\infty}[v_j(s)-u_{\theta}]| < \epsilon,
\, s \in p(2 \delta;t), \, j \in \{ 1,2,\ldots, N \}, \, k > K.
\]
Hence,
\begin{align*}
& \left\| \int_0^t \mathbf{\omega}  \{ S_{\beta_k}[\mathbf{u}_{\beta_k}(s)-\mathbf{u_{\theta}}]
- S_{\infty}[\mathbf{v}(s)-\mathbf{u_{\theta}}] \} \, ds \right\|_{\infty} \\
& \mbox{ } = \left\| \int_{p(2\delta;t) \, \cup \, r(2\delta;t)} \mathbf{\omega}  \{ S_{\beta_k}[\mathbf{u}_{\beta_k}(s)-\mathbf{u_{\theta}}]
- S_{\infty}[\mathbf{v}(s)-\mathbf{u_{\theta}}] \} \, ds \right\|_{\infty} \\
& \mbox{ } \leq \left\| \int_{p(2\delta;t)} \mathbf{\omega}  \{ S_{\beta_k}[\mathbf{u}_{\beta_k}(s)-\mathbf{u_{\theta}}]
- S_{\infty}[\mathbf{v}(s)-\mathbf{u_{\theta}}] \} \, ds \right\|_{\infty} \\
& \mbox{ } \hspace*{0.5cm}+\left\| \int_{r(2\delta;t)} \mathbf{\omega}  \{ S_{\beta_k}[\mathbf{u}_{\beta_k}(s)-\mathbf{u_{\theta}}]
- S_{\infty}[\mathbf{v}(s)-\mathbf{u_{\theta}}] \} \, ds \right\|_{\infty} \\
& \mbox{ } \leq \epsilon |p(2\delta;t)| \max_{i \in \{1,2,\ldots,N\}} \left( \sum_{j=1}^{N}  |\omega_{i,j}| \right) \\
& \mbox{ } \hspace*{0.5cm}+ |r(2\delta;t)| \max_{i \in \{1,2,\ldots,N\}} \left( \sum_{j=1}^{N}  |\omega_{i,j}| \right) \\
& \mbox{ } \leq \frac{\tilde{\epsilon}}{2TC} T \max_{i \in \{1,2,\ldots,N\}} \left( \sum_{j=1}^{N}  |\omega_{i,j}| \right) \\
& \mbox{ } \hspace*{0.5cm}+\frac{\tilde{\epsilon}}{2C}  \max_{i \in \{1,2,\ldots,N\}} \left( \sum_{j=1}^{N}  |\omega_{i,j}| \right) \\
& \mbox{ } < \tilde{\epsilon} 
\end{align*}
for all $k > K$, where the second last inequality follows from (\ref{C0.01}), the fact that $|p(2\delta;t)| \leq T$ for $t \in [0,T]$ and (\ref{C0}). Since $\tilde{\epsilon} > 0$ and $t \in [0,T]$ were arbitrary, we conclude that (\ref{C5}) must hold. \\
\rule{3mm}{3mm}


\subsection{Limit problem}
By employing the uniform convergence (\ref{limit-function}), the convergence of the integral (\ref{C5}) and assumption (\ref{limit-int-q}), we conclude from (\ref{integral-form}) that the limit function $\mathbf{v}$ satisfies
\begin{align}
\nonumber
\mathbf{\tau} \mathbf{v}(t) - \mathbf{\tau} \mathbf{u}_{\mathrm{init}} = & - \int_0^t \mathbf{v}(s) \, ds \\
\nonumber
& + \int_0^t \mathbf{\omega}  S_{\infty}[\mathbf{v}(s)-\mathbf{u_{\theta}}] \, ds \\
\label{C5.1}
& + \int_0^t \mathbf{q}_{\infty}(s) \, ds, \, t \in [0,T],
\end{align}
provided that $\mathbf{v}$ is threshold simple.
Recall that $\mathbf{v}$ is continuous. Consequently, if $\mathbf{v}$ is extra threshold simple, then it follows from the Fundamental Theorem of Calculus that $\mathbf{v}$ also satisfies the ODEs, except at time instances where one or more of the component functions equals the threshold value for firing: 
\begin{align}
\label{C6}
\mathbf{\tau} \mathbf{v}'(t) =& - \mathbf{v}(t) + \mathbf{\omega}  S_{\infty}[\mathbf{v}(t)-\mathbf{u_{\theta}}] + \mathbf{q}_{\infty}(t), \, t \in (0,T] \setminus Z(\mathbf{v}), \\
\label{C7}
\mathbf{v}(0)=&\mathbf{u}_{\mathrm{init}},
\end{align}
where $Z(\mathbf{v})$ is defined in (\ref{zero-measure}). 

The existence of a solution matter, for point neuron models with a Heaviside firing rate function, is summarized in the following theorem: 
\begin{theorem} 
\label{theorem-existence}
If the limit $\mathbf{v}$ in (\ref{limit-function}) is threshold simple, then $\mathbf{v}$ solves (\ref{C5.1}). 
In the case that $\mathbf{v}$ is extra threshold simple, $\mathbf{v}$  also satisfies (\ref{C6})-(\ref{C7}).
\end{theorem}
In \cite{PG2010} the existence issue for neural field equations with a Heaviside activation function is studied, but the analysis is different because a continuum model is considered. We would also like to mention that Theorem \ref{theorem-existence} can not be regarded as a simple consequence of Carath\'{e}odory's existence theorem \cite{WikipediaCaratheodory} because the right-hand-side of (\ref{C6}) is discontinuous with respect to $\mathbf{v}$.

\section{Uniqueness}
\label{uniqueness}
If $\beta < \infty$, then standard ODE theory \cite{WikipediaPicard} implies that (\ref{point-neuron-model})-(\ref{initial-condition}) has a unique solution. Unfortunately, as will be demonstrated below, this desirable property is not necessarily inherited by the infinite steepness limit problem. 

We will first explain why the uniqueness question is a subtle issue for point neuron models with a Heaviside firing rate function. Thereafter, additional requirements are introduced which ensure the uniqueness of an extra threshold simple solution.  

\subsection{Example: Several solutions}
\label{example-several-solutions}
Let us study the problem
\begin{align}
\label{D1}
v'(t) &= -v(t) + \omega S_{\infty}[v(t)-u_{\theta}], \quad t \in (0,T], \\
\label{D2}
v(0) &= u_{\theta},
\end{align}
where we assume that \[ w>u_{\theta} \geq 0.\]
Note that the ODE (\ref{D1}) is not required to hold for $t=0$.
Consider the functions
\begin{align}
\label{v1}
v_1(t) &= \omega + (u_{\theta}-\omega) e^{-t} = u_{\theta} e^{-t} +(1-e^{-t}) \omega, \\
\label{v2}
v_2(t) &= u_{\theta} e^{-t}.
\end{align}
Since
\begin{eqnarray*}
v_1(t) &>& u_{\theta} e^{-t} +(1-e^{-t}) u_{\theta}  = u_{\theta}, \quad t \in (0,T], \\
v_2(t) &<& u_{\theta}, \quad t \in (0,T],
\end{eqnarray*}
it follows that both $v_1$ and $v_2$ solves (\ref{D1})-(\ref{D2}). 

Furthermore, with 
\begin{align*} 
&\omega=2 u_{\theta}, \\
&S_{\infty}(0)=\frac{1}{2}, 
\end{align*}
we actually obtain a third solution of (\ref{D1})-(\ref{D2}). 
More specifically, the stationary solution 
\begin{equation}
\label{v3}
v_3(t)=u_{\theta}, \, t \in [0,T].
\end{equation}

We conclude that models with a Heaviside firing rate function can have several solutions -- such problems can thus become ill-posed.
(In \cite{NW2016} we showed that the {\em initial-condition-to-solution} map is not necessarily continuous for such problems, and that the error amplification ratio can become very large in the steep, but Lipschitz continuous, firing rate regime).
Note that switching to the integral form (\ref{C5.1}) will not resolve the lack of uniqueness issue for the toy example considered in this subsection.

We also remark that: 
\begin{itemize}
\item If we define $S_{\infty}(0)=1/2$, then neither $v_1$ nor $v_2$ satisfies the ODE (\ref{D1}) for $t=0$. (In the case $\omega=2 u_{\theta}$, $v_3$ satisfies (\ref{D1}) for $t=0$.)
\item If we define $S_{\infty}(0)=1$, then $v_1$, but not $v_2$, satisfies (\ref{D1}) also for $t=0$.
\item If we define $S_{\infty}(0)=0$, then $v_2$, but not $v_1$, satisfies (\ref{D1}) also for $t=0$.
\end{itemize}

\subsection{Enforcing uniqueness}
In order to enforce uniqueness, we need to impose further restrictions.
It turns out that it is sufficient to require that the derivative is continuous from the right and that the ODEs also must be satisfied whenever one, or more, of the component functions equals the threshold value for firing:
\begin{align}
\label{E1}
\mathbf{\tau} \mathbf{v}'(t) &= - \mathbf{v}(t) + \mathbf{\omega}  S_{\infty}[\mathbf{v}(t)-\mathbf{u_{\theta}}] + \mathbf{q}_{\infty}(t), \, t \in [0,T],\\
\label{E2}
\mathbf{v}(0)&=\mathbf{u}_{\mathrm{init}}.
\end{align}
Note that the ODEs (\ref{E1}) also must be satisfied for $t=0$, in case one of the components of $\mathbf{u}_{\mathrm{init}}$ equals $u_{\theta}$.

\vspace*{0.2cm}
\noindent{\bf Definition} {\bf \em (Right smooth)} \\
A vector-valued function $\mathbf{z}:[0,T] \rightarrow \Re^N$ is right smooth if  $\mathbf{z}'$ is continuous 
from the right for all $t \in [0,T)$. 
\vspace*{0.2cm}

\begin{theorem}
\label{theorem-uniqueness}
Equations (\ref{E1})-(\ref{E2}) can at the most have one solution which is both extra threshold simple and right smooth. 
\end{theorem}

\subsubsection*{Proof}
Let $\mathbf{v}$ and $\mathbf{\tilde{v}}$ be two solutions of (\ref{E1})-(\ref{E2}) which are both right smooth and extra threshold simple: 
\begin{eqnarray*}
&& [0,T] = \bigcup_{l=1}^{L} \bar{I}_l, \\
&& m(s;\mathbf{v}) \neq 0 \quad \forall s \in \bigcup_{l=1}^{L}  I_l,
\end{eqnarray*}
and 
\begin{eqnarray*}
&& [0,T] = \bigcup_{l=1}^{\tilde{L}} \bar{\tilde{I}}_l, \\
&& m(s;\mathbf{\tilde{v}}) \neq 0 \quad \forall s \in \bigcup_{l=1}^{\tilde{L}} \tilde{I}_l.
\end{eqnarray*}
where $I_1, I_2, \ldots, I_{L}$ and $\tilde{I}_1, \tilde{I}_2, \ldots, \tilde{I}_{\tilde{L}}$ are disjoint open intervals,
see (\ref{ms}) and the definition of extra threshold simple in subsection \ref{threshold-terminology}.

Then there exist disjoint open intervals $\hat{I}_1, \hat{I}_2, \ldots, \hat{I}_{\hat{L}}$
such that
\begin{align}
\nonumber
& [0,T] = \bigcup_{l=1}^{\hat{L}} \bar{\hat{I}}_l, \\
\label{E2.5}
& m(s;\mathbf{v}) \neq 0 \mbox{ and } m(s;\mathbf{\tilde{v}}) \neq 0 \quad \forall s \in \bigcup_{l=1}^{\hat{L}} \hat{I}_l.
\end{align}
Let us focus on one of these intervals, $\hat{I}_l = (a_l,a_{l+1})$.
Define 
\[ 
\mathbf{d}=\mathbf{v}-\mathbf{\tilde{v}} 
\] 
and assume that
\begin{equation}
\label{E2.7}
\mathbf{v}(a_l) = \mathbf{\tilde{v}}(a_l), 
\end{equation}
which obviously holds for $l=1$.
Then,
\begin{align}
\label{E3}
\mathbf{\tau} \mathbf{d}'(t) &= - \mathbf{d}(t) + \mathbf{\omega} \mathbf{\gamma}(t),
\, t \in [a_l,a_{l+1}],\\
\label{E4}
\mathbf{d}(a_l)&=\mathbf{0},
\end{align}
where
\[
\mathbf{\gamma}(t) = S_{\infty}[\mathbf{v}(t)-\mathbf{u_{\theta}}] - S_{\infty}[\mathbf{\tilde{v}}(t)-\mathbf{u_{\theta}}], \, t \in [a_l,a_{l+1}].
\]
Note that, due to (\ref{E2.5}), $\mathbf{\gamma}(t)$ equals a constant vector $\mathbf{c}$, with components $-1,0$ or $1$, except possibly at $t=a_l, \, a_{l+1}$:
\begin{equation}
\label{E4.1}
\mathbf{\gamma}(t) = \mathbf{c}, \, t \in (a_l,a_{l+1}).
\end{equation}
Furthermore, from (\ref{E2.7}) we find that 
\begin{equation}
\label{E5}
\mathbf{\gamma}(a_l) = \mathbf{0}.
\end{equation}

Putting $t=a_l$ in (\ref{E3}) and invoking (\ref{E4}) and (\ref{E5}) yield that
\[
\mathbf{d}'(a_l) = \mathbf{0},
\]
and from the right continuity of $\mathbf{d}'$ and $\mathbf{d}$, (\ref{E3}), (\ref{E4}) and (\ref{E4.1}) we find that
\[
\mathbf{0} = \mathbf{\tau} \mathbf{d}'(a_l)
= \lim_{t \rightarrow a_l^+} \mathbf{\tau} \mathbf{d}'(t)
= \lim_{t \rightarrow a_l^+} [- \mathbf{d}(t) + \mathbf{\omega} \mathbf{\gamma}(t)]
= \mathbf{\omega} \mathbf{c}.
\]
Since $\mathbf{\omega} \mathbf{\gamma}(t) = \mathbf{\omega} \mathbf{c} = \mathbf{0}$, $t \in (a_l, a_{l+1})$,  and $\mathbf{\omega} \mathbf{\gamma}(a_l) =\mathbf{0} $, see (\ref{E5}), we conclude from (\ref{E3})-(\ref{E4}) that $\mathbf{d}$ satisfies 
\begin{eqnarray*}
\mathbf{\tau} \mathbf{d}'(t) &=& - \mathbf{d}(t),
\, t \in [a_l,a_{l+1}),\\
\mathbf{d}(a_l)&=&\mathbf{0},
\end{eqnarray*}
which has the unique solution $\mathbf{d}(t) = \mathbf{0}$, $t \in [a_l,a_{l+1})$. 
Both $\mathbf{v}(t)$ and $\mathbf{\tilde{v}}(t)$ are differentiable on $[0,T]$ and hence continuous.  
It follows that, by employing the continuity of $\mathbf{v}$ and $\mathbf{\tilde{v}}$ at time $t=a_{l+1}$,
\[
\mathbf{v}(t) = \mathbf{\tilde{v}} (t), \, t \in [a_l,a_{l+1}].
\]

Since $\mathbf{v}(a_{l+1}) = \mathbf{\tilde{v}} (a_{l+1})$, we can repeat the argument on the next interval
$[a_{l+1},a_{l+2}]$, and it follows by induction that $\mathbf{v}(t) = \mathbf{\tilde{v}} (t), \, t \in [0,T]$. 
\\ \rule{3mm}{3mm} 

We would like to comment the findings presented in the bullet-points at the end of subsection \ref{example-several-solutions} in view of Theorem \ref{theorem-uniqueness}: In order to enforce uniqueness for the solution of 
(\ref{D1})-(\ref{D2}), we can require that the ODE (\ref{D1}) also should be satisfied for $t=0$. Nevertheless, this might force us to define $S_{\infty}(0) \neq \frac{1}{2}$,  which differs from the standard definition of the Heaviside function $H$. 

More generally, if one has accomplished to compute an extra threshold simple and right smooth function $\mathbf{v}$ which satisfies (\ref{C5.1}), then one can attempt to redefine $S_{\infty}[\mathbf{v}(t)-\mathbf{u_{\theta}}]$, $t \in \{a_1, \, a_2, \ldots, a_{L+1} \}$, 
such that (\ref{E1})-(\ref{E2}) hold, and $\mathbf{v}$ is the only solution to this problem. This may imply that $S_{\infty}[\mathbf{v}(t)-\mathbf{u_{\theta}}]$ can not be generated by using the composition $H \circ [\mathbf{v}(t)-\mathbf{u_{\theta}}]$. Instead one must determine $z_{j,k}=S_{\infty} [v_j(a_k)-u_{\theta}]$, $j=1, \, 2, \ldots, \, N$, $k = 1, \, 2, \ldots, \, L+1$. More precisely, for each $k \in \{1, \, 2, \ldots, \, L+1 \}$ one gets a  linear system of algebraic equations
\[
\tau_{i} v'_i(a_k) = -v_i(a_k) + \sum_{j=1}^{N} \omega_{i,j} z_{j,k} + q_{\infty,i}(a_k), \, i=1,2,\ldots,N, 
\]
which will have a unique solution $(z_{1,k}, \, z_{2,k}, \, \ldots, \, z_{N,k})^T$ if the connectivity matrix $\mathbf{\omega} = [\omega_{i,j}]$ is nonsingular. (In this paragraph, $\{0=a_1, \, a_2, \ldots, a_{L+1}=T \}$ are the time instances one, or more, of the component functions of $\mathbf{v}$ potentially equals the threshold value for firing, see the definition of extra threshold simple in subsection \ref{threshold-terminology}).   

\section{Convergence of the entire sequence}
\label{convergence-entire-sequence} 
We have seen that point neuron models with a Heaviside firing rate function can have several solutions. 
One therefore might wonder, can different subsequences of $\{ \mathbf{u}_{\beta} \}$ converge to different solutions of the limit problem? In this section we present an example which shows that this can happen, even though the involved sigmoid functions  
satisfy assumption {\bf \cal{A}}. 

\subsection{Example: Different subsequences can converge to different solutions}
Let us again consider the initial value problem (\ref{D1})-(\ref{D2}), which we discussed in subsection \ref{example-several-solutions}. 
A finite steepness approximation of this problem, using the notation $u(t) =u_{\beta} (t)$, reads: 
\begin{align}
\label{G1}
u'(t) &= -u(t) + \omega \bar{S}_{\beta}[u(t)-u_{\theta}], \quad t \in (0,T], \\
\label{G2}
u(0) &= u_{\theta},
\end{align}
where 
\[
\bar{S}_{\beta}[x] = S_{\beta} \left[ x+\frac{(-1)^{\beta}}{2\beta} \right], \, \beta \in \naturalNumbers,  
\]
and $S_{\beta}$ is, e.g., either the $\tanh$-based sigmoid function (\ref{firingtanh1})-(\ref{firingtanh2}) or 
\begin{equation}
\label{alt-firing}
S_{\beta} (x)= \left\{
\begin{array}{ll}
1,& x > \frac{1}{\beta}, \\
\frac{1}{2}+\frac{1}{2} \beta x ,& x \in \left[ -\frac{1}{\beta}, \frac{1}{\beta} \right], \\
0, & x < -\frac{1}{\beta}.
\end{array}
\right.
\end{equation} 
Note that $\{ \bar{S}_{\beta} \}$ converges pointwise, except for $x=0$, to the Heaviside function $H$ as $\beta \rightarrow \infty$. 
In fact, $\{ \bar{S}_{\beta} \}$ satisfies assumption {\bf \cal{A}}.  

\begin{figure}
  \centering
  \begin{subfigure}[b]{\textwidth}
  \includegraphics[width=15cm,trim=20mm 95mm -20mm 90mm,clip=true]{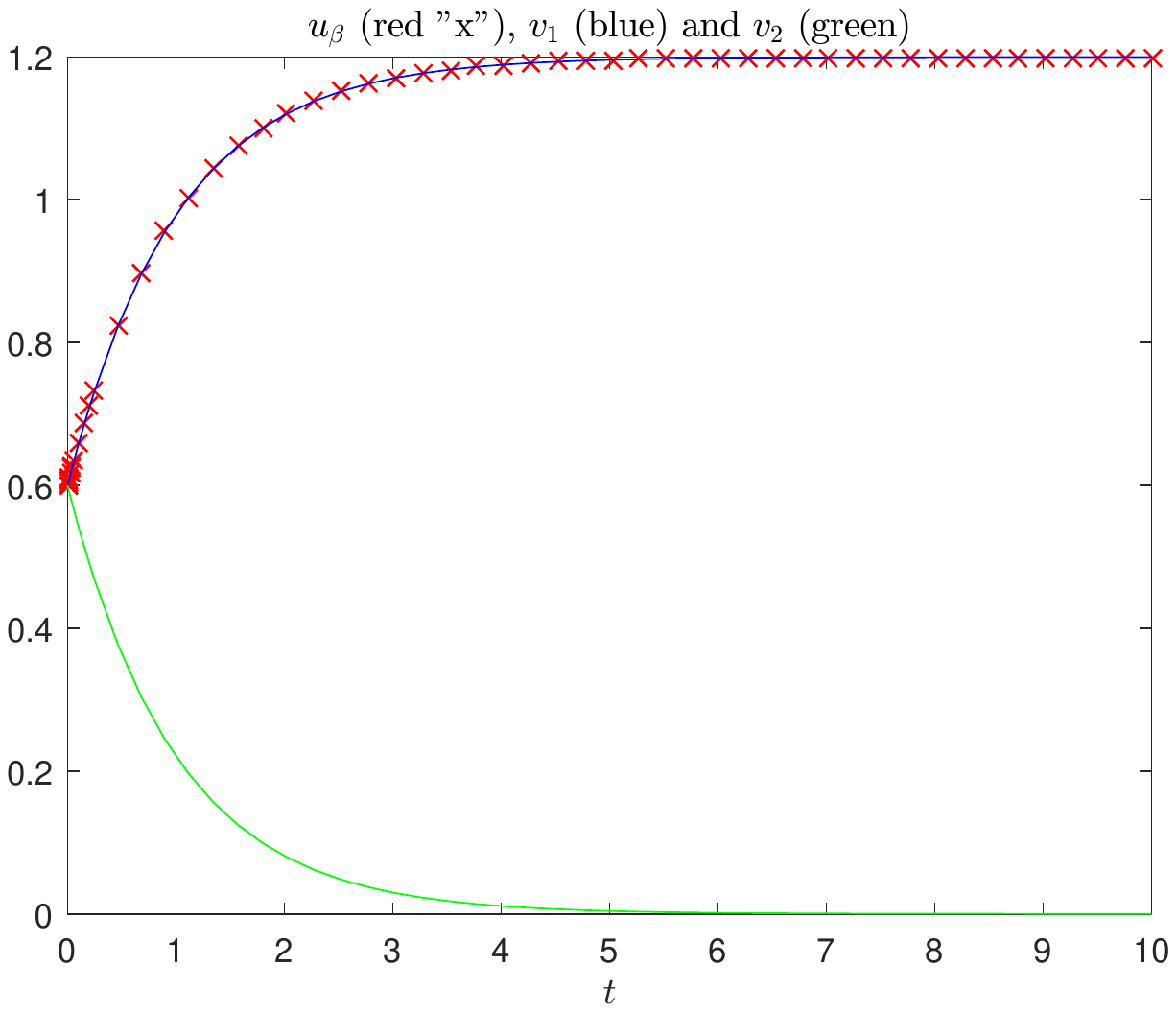}
  \caption{$\beta=10\, 000\, 000$}
  \end{subfigure}

  \begin{subfigure}[b]{\textwidth}
  \includegraphics[width=15cm,trim=20mm 95mm -20mm 90mm,clip=true]{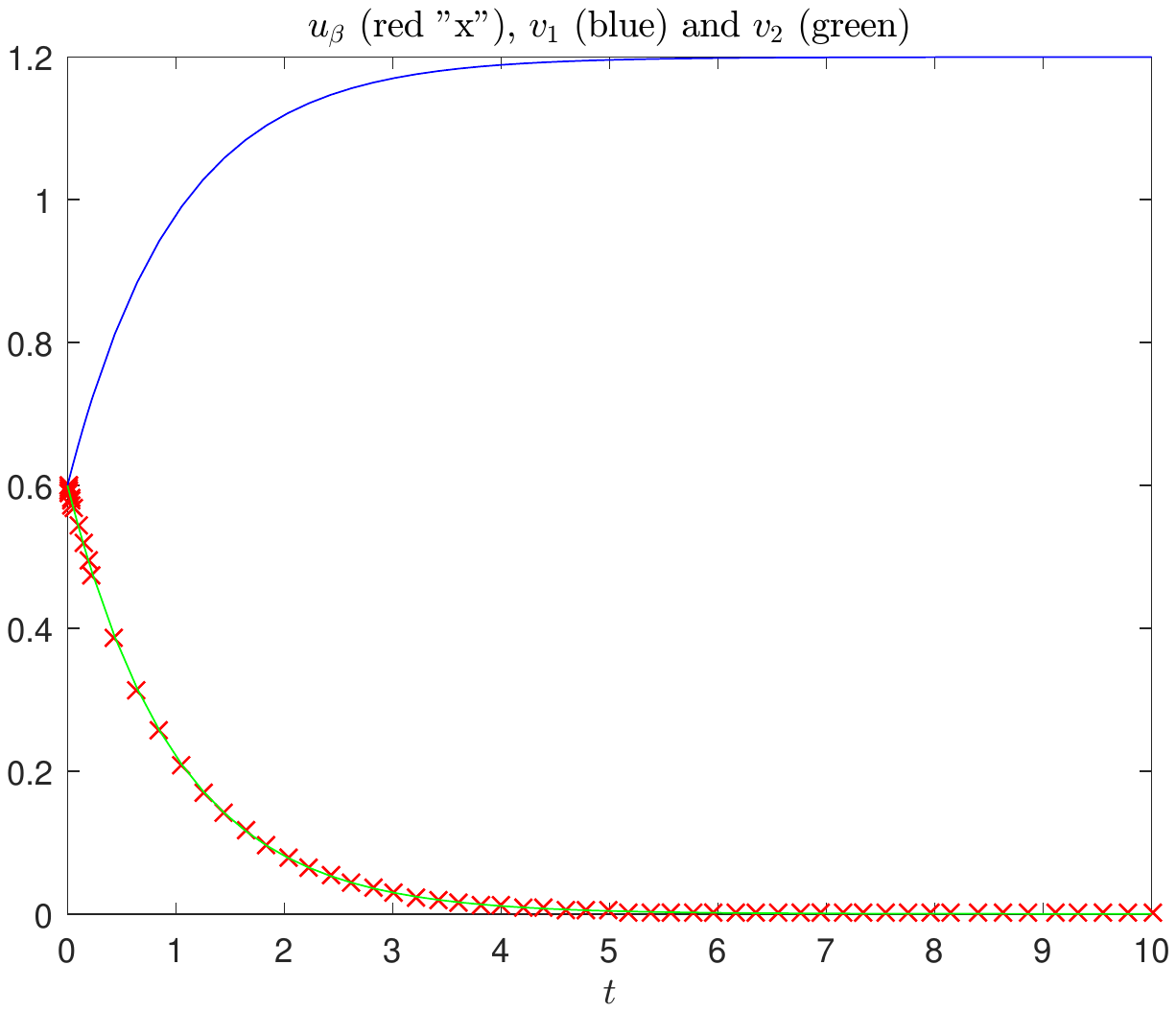}
  \caption{$\beta=10\, 000\, 001$}
  \end{subfigure}
  \caption{Numerical solutions of (\ref{G1})-(\ref{G2}) computed with Matlab's {\tt ode45} software. In these simulations we used   $u_{\theta}=0.6$ and $\omega=1.2$. The functions $v_1$ and $v_2$, see (\ref{v1}) and (\ref{v2}), are the solutions of the associated limit problem (\ref{D1})-(\ref{D2}).}\label{figure-numerical}
\end{figure}

We consider the case $\omega=2 u_{\theta}$, and (\ref{D1})-(\ref{D2}) therefore has three solutions $v_1, \, v_2$ and $v_3$, see (\ref{v1}), (\ref{v2}) and (\ref{v3}) in subsection \ref{example-several-solutions}. Note that \[ u(t) =u_{\beta} (t) \] has the property 
\begin{itemize}
\item $u'_{\beta}(0) > c$ if $\beta$ is even, 
\item $u'_{\beta}(0) < -c$ if $\beta$ is odd,  
\end{itemize}
where $c>0$ is a constant which is independent of $\beta$. 
It therefore follows, argument not included, that 
\begin{align}
& \lim_{k \rightarrow \infty} u_{2k} = v_1, \\
& \lim_{k \rightarrow \infty} u_{2k+1} = v_2, 
\end{align}
and no subsequence converges to the third solution $v_3$. Figure \ref{figure-numerical} shows numerical solutions of 
(\ref{G1})-(\ref{G2}) with steepness parameter $\beta = 10\, 000\, 000, \\ \, 10\, 000\, 001$, using the firing rate function (\ref{alt-firing}) to define $\bar{S}_{\beta}$. (If one instead employs (\ref{firingtanh1})-(\ref{firingtanh2}) in the implementation of $\bar{S}_{\beta}$, the plots, which are not included, are virtually unchanged). 

We would like to mention that we have not been able to construct an example of this kind for Lipschitz continuous firing rate functions which converge pointwise to the Heaviside function also for $x=0$. 

\subsection{Entire sequence}
We have seen that almost everywhere convergence of the sequence of firing rate functions to the Heaviside limit is not sufficient to guarantee that the entire sequence $\{ u_{\beta} \}$ converges to the same solution of the limit problem. Nevertheless, one has the following result: 
\begin{theorem}
\label{theorem-entire-convergence}
Let $\mathbf{v}$ be the limit function in (\ref{limit-function}). 
If the limit of every convergent subsequence of $\{ \mathbf{u}_{\beta} \}$ is extra threshold simple, right smooth and 
satisfies  (\ref{E1})-(\ref{E2}), 
then the entire sequence $\{ \mathbf{u}_{\beta} \}$ converges uniformly to $\mathbf{v}$.
\end{theorem}
\subsubsection*{Proof}
Suppose that the entire sequence $\{ \mathbf{u}_{\beta} \}$ does  \underline{not} converge uniformly to $\mathbf{v}$. 
Then there is an $\epsilon > 0$ such that, for every positive integer $M$, there must exist $\mathbf{u}_{\beta_l}$, $\beta_l > M$, satisfying  
\begin{equation}
\label{F1}
\sup_{t \in [0,T]} \|\mathbf{u}_{\beta_l}(t) -\mathbf{v}(t) \|_{\infty} > \epsilon.
\end{equation}
The subsequence $\{ \mathbf{u}_{\beta_l} \}$ can thus not converge uniformly to $\mathbf{v}$, but constitute a set of uniformly bounded and equicontinuous functions, see section \ref{uniformily-bounded-equicontinuous}. According to the Arzel\`{a}-Ascoli theorem, $\{ \mathbf{u}_{\beta_l} \}$ therefore possesses a  uniformly convergent subsequence 
 $\{ \mathbf{u}_{\beta_{l_n}} \}$, 
 \[
 \lim_{n \rightarrow \infty} \mathbf{u}_{\beta_{l_n}} = \mathbf{\tilde{v}}.
 \]
 Due to (\ref{F1}), 
 \begin{equation}
 \label{F2}
 \mathbf{\tilde{v}} \neq \mathbf{v}. 
 \end{equation}
 
On the other hand, both $\mathbf{v}$ and $ \mathbf{\tilde{v}}$ are limits of subsequences of $\{ \mathbf{u}_{\beta} \}$, and are therefore by assumption  extra threshold simple, right smooth and satisfies  (\ref{E1})-(\ref{E2}). Hence, Theorem \ref{theorem-uniqueness} implies that 
 $ \mathbf{\tilde{v}} = \mathbf{v}$, which contradicts (\ref{F2}).  
We conclude that the entire sequence $\{ \mathbf{u}_{\beta} \}$ must converge uniformly to $\mathbf{v}$.
\\ \rule{3mm}{3mm}

One might argue that Theorem \ref{theorem-entire-convergence} only is of theoretical interest because it seems very difficult to guarantee that ``the limit of every convergent subsequence of $\{ \mathbf{u}_{\beta} \}$ is extra threshold simple, right smooth and satisfies  (\ref{E1})-(\ref{E2})''. 

\section{Example: Threshold advanced limits} 
\label{threshold-advanced-limits}
We will now show that threshold advanced limits, i.e. limits which are not threshold simple, may possess some peculiar properties. More precisely, such limits can potentially occur in (\ref{limit-function}), and they do not necessarily satisfy the limit problem obtained by using a Heaviside firing rate function. 

With source terms which do not depend on the steepness parameter $\beta$, we have not managed to construct an example with a threshold advanced limit $\mathbf{v}$. If we allow $\mathbf{q}=\mathbf{q}_{\beta}$, this can, however, be accomplished as follows. 
Let 
\[
z_{\beta}(t) = \frac{1}{\beta} S_{\beta}[-\frac{1}{\beta} + 2t] + u_{\theta}, 
\, \beta = 1,2,\ldots, 
\]
where we, for the sake of simplicity, work with the firing rate function (\ref{alt-firing}). 
Then, 
\begin{eqnarray*}
&& z_{\beta}(0) = \frac{1}{\beta} S_{\beta}[-\frac{1}{\beta}] + u_{\theta} = u_{\theta}, \\
&& z_{\beta}(t)  = 
\left\{ 
\begin{array}{ll}
t + u_{\theta}, & t \in [0,\frac{1}{\beta})  \\
\frac{1}{\beta} +  u_{\theta}, & t \geq \frac{1}{\beta}, 
\end{array}
\right. \\
&& z'_{\beta}(t)  = 
\left\{ 
\begin{array}{ll}
1, & t \in [0,\frac{1}{\beta})  \\
0, & t > \frac{1}{\beta}, 
\end{array}
\right. \\
&& S_{\beta}[z_{\beta}(t)- u_{\theta}]=
\left\{ 
\begin{array}{ll}
\frac{1}{2}+\frac{1}{2} \beta t, & t \in [0,\frac{1}{\beta}) \\
1, & t \geq \frac{1}{\beta}, 
\end{array}
\right. \\
\end{eqnarray*}
and we find that 
\[
u_{\beta}(t) = z_{\beta}(t)
\] 
solves 
\begin{eqnarray*}
u_{\beta}(t) - u_{\theta}&=& - \int_0^t u_{\beta}(s) \, ds \\
&& + \int_0^t \omega  S_{\beta}[u_{\beta}(s)-u_{\theta}] \, ds \\
&& + \int_0^t q_{\beta}(s) \, ds, \, t \in [0,T],
\end{eqnarray*}
where 
\begin{align}
\nonumber
q_{\beta}(t) &= z'_{\beta}(t)+z_{\beta}(t)-\omega S_{\beta}[z_{\beta}(s)-u_{\theta}] \\
&= 
\label{discontinuous-source}
\left\{
\begin{array}{ll}
1+t+u_{\theta}-\omega (\frac{1}{2}+\frac{1}{2} \beta t), & t \in [0,\frac{1}{\beta}) \\
\frac{1}{\beta}+u_{\theta}-\omega, & t > \frac{1}{\beta}. 
\end{array}
\right.
\end{align}

It follows that 
\[
q_{\infty}(t)  = 
\left\{
\begin{array}{ll}
1+u_{\theta} - \omega, & t=0 \\
u_{\theta} - \omega, & t>0,
\end{array}
\right.
\]
and since, for any $\beta \in \naturalNumbers$,   
\[
|q_{\beta}(t) | \leq 1+\frac{1}{\beta}+|u_{\theta}|+|\omega|
< 2 + |u_{\theta}|+|\omega|, \, t \neq \frac{1}{\beta}, 
\]
we conclude that   
\[
\lim_{\beta \rightarrow \infty} \int_0^t q_{\beta}(s) \, ds = \int_0^t q_{\infty}(s) \, ds , \, t \in [0,T].
\]

Note that 
\[
u_{\beta} (t) \longrightarrow \bar{v}(t)=u_{\theta}, \mbox{ uniformly, as }  \beta \rightarrow \infty, 
\]
but $\bar{v}(t)=u_{\theta}$ does not solve the limit problem 
\begin{eqnarray*}
v(t) - u_{\theta}&=& - \int_0^t v(s) \, ds \\
&& + \int_0^t \omega  S_{\infty}[v(s)-u_{\theta}] \, ds \\
&& + \int_0^t q_{\infty}(s) \, ds, \, t \in [0,T],
\end{eqnarray*}
because 
\begin{eqnarray*}
&& - \int_0^t \bar{v}(s) \, ds 
+ \int_0^t \omega  S_{\infty}[\bar{v}(s)-u_{\theta}] \, ds 
+ \int_0^t q_{\infty}(s) \, ds \\
&& = -t u_{\theta}+t \omega \frac{1}{2} + t (u_{\theta} - \omega) \\
&& = - \frac{1}{2} t \omega \\
&& \neq 0 = \bar{v}(t) - u_{\theta} , \, t \in (0,T].
\end{eqnarray*}
This argument assumes that $S_{\infty}[0] = 1/2$. If one instead defines $S_{\infty}[0] = 1$, then $\bar{v}$ would solve 
the limit problem. 

Due to the properties of the firing rate function (\ref{alt-firing}), the source term $q_{\beta}$ in (\ref{discontinuous-source}) becomes discontinuous. This can be avoided by instead using the smooth version (\ref{firingtanh1})-(\ref{firingtanh2}), but then the analysis of this example becomes much more involved.  

The author does not know whether it is possible to impose restrictions which would guarantee that the limit $\mathbf{v}$ in 
(\ref{limit-function}) is threshold simple or extra threshold simple. This seems to be an herculean task. 

\section{Discussion and conclusions}
If a Heaviside firing rate function is used, then the model (\ref{model-eq1})-(\ref{model-eq2}) may not only have several solutions, but the {\em initial-condition}-to-{\em solution} map for this problem can become discontinuous \cite{NW2016}. It is thus virtually impossible to develop reliable numerical methods, which employ finite precision arithmetic, for such problems. One can try to overcome this issue by: 
\begin{description}
\item[a)] Attempting to solve the ill-posed equation with symbolic computations. 
\item[b)] Regularize the problem. 
\end{description}
As far as the author knows, present symbolic techniques are not able to handle  strongly nonlinear equations of the kind (\ref{model-eq1}), even when $\beta < \infty$. We therefore analyzed the approach b), using the straightforward regularization technique obtained by replacing the Heaviside firing rate function by a Lipschitz continuous mapping. This  yields an equation which is within the scope of the Picard-Lindel\"{o}f theorem and standard stability estimates for ODEs. That is, well-posed and, at least in principle,  approximately solvable by numerical methods. 

Our results show that the sequence $\{ \mathbf{u}_{\beta} \}$ of regularized solutions will have at least one convergent subsequence. The limit, $\mathbf{v}$, of this subsequence will satisfy the integral/Volterra form (\ref{C5.1}) of the limit problem, provided that the Lebesgue measure of the time one, or more, of the component functions of $\mathbf{v}$ equals the threshold value $u_{\theta}$ for firing, is zero. Unfortunately, it seems to be very difficult to impose restrictions which would guarantee that $\mathbf{v}$ obeys this threshold property, which we refer to as threshold simple. Also, the example presented in section \ref{threshold-advanced-limits} shows that, if the limit $\mathbf{v}$ is not threshold simple, then this function may not solve the associated equation with a Heaviside firing rate function. 

One could propose to overcome the difficulties arising when $\beta = \infty$ by always working with finite slope firing rate functions. This would potentially yield a rather robust approach, provided that the entire sequence $\{ \mathbf{u}_{\beta} \}$ converges, because increasing a large $\beta$ would still guarantee that $\mathbf{u}_{\beta}$ is close to the unique limit $\mathbf{v}$.  
However, the fact that different convergent subsequences of $\{ \mathbf{u}_{\beta} \}$ can converge to different solutions of the limit problem, as discussed in section \ref{convergence-entire-sequence}, suggests that this approach must be applied with great care. In addition, the error amplification in the steep firing rate regime can become extreme \cite{NW2016}, and the accurate numerical solution of such models is thus challenging.  

What are the practical consequences of our findings? 
As long as there does not exist very reliable biological information about the size of the steepness parameter $\beta$, and the shape of the firing rate function $S_{\beta}$, it seems that we have to be content with simulating with various $\beta < \infty$. If one observes that $\mathbf{u}_{\beta}$ approaches a threshold advanced limit, as $\beta$ increases, or that the entire sequence does not converge, the alarm bell should ring. All simulations with large $\beta$ must use error control methods which guarantee the accuracy of the numerical solution --  we must keep in mind that we are trying to solve an almost ill-posed problem.

\section*{Competing interests}
The author declares that he has no competing interests.

\section*{Acknowledgements}
This work was supported by the The Research Council of Norway, project number 239070. 
The author would like to thank Prof. Wyller for several interesting discussions about the research presented in this paper. 

\end{document}